\documentclass[10pt]{article}
\usepackage{times,amsmath,amssymb,theorem}

\newdimen\proofrulebreadth \proofrulebreadth=.05em
\newdimen\proofdotseparation \proofdotseparation=1.25ex
\newdimen\proofrulebaseline \proofrulebaseline=2ex
\newcount\proofdotnumber \proofdotnumber=3
\let\then\relax
\def\hfi{\hskip0pt plus.0001fil}
\mathchardef\squigto="3A3B
\newif\ifinsideprooftree\insideprooftreefalse
\newif\ifonleftofproofrule\onleftofproofrulefalse
\newif\ifproofdots\proofdotsfalse
\newif\ifdoubleproof\doubleprooffalse
\let\wereinproofbit\relax
\newdimen\shortenproofleft
\newdimen\shortenproofright
\newdimen\proofbelowshift
\newbox\proofabove
\newbox\proofbelow
\newbox\proofrulename
\def\shiftproofbelow{\let\next\relax\afterassignment\setshiftproofbelow\dimen0 }
\def\shiftproofbelowneg{\def\next{\multiply\dimen0 by-1 }%
\afterassignment\setshiftproofbelow\dimen0 }
\def\setshiftproofbelow{\next\proofbelowshift=\dimen0 }
\def\setproofrulebreadth{\proofrulebreadth}

\def\prooftree{
\ifnum  \lastpenalty=1
\then   \unpenalty
\else   \onleftofproofrulefalse
\fi
\ifonleftofproofrule
\else   \ifinsideprooftree
        \then   \hskip.5em plus1fil
        \fi
\fi
\bgroup
\setbox\proofbelow=\hbox{}\setbox\proofrulename=\hbox{}%
\let\justifies\proofover\let\leadsto\proofoverdots\let\Justifies\proofoverdbl
\let\using\proofusing\let\[\prooftree
\ifinsideprooftree\let\]\endprooftree\fi
\proofdotsfalse\doubleprooffalse
\let\thickness\setproofrulebreadth
\let\shiftright\shiftproofbelow \let\shift\shiftproofbelow
\let\shiftleft\shiftproofbelowneg
\let\ifwasinsideprooftree\ifinsideprooftree
\insideprooftreetrue
\setbox\proofabove=\hbox\bgroup$\displaystyle 
\let\wereinproofbit\prooftree
\shortenproofleft=0pt \shortenproofright=0pt \proofbelowshift=0pt
\onleftofproofruletrue\penalty1
}

\def\eproofbit{
\ifx    \wereinproofbit\prooftree
\then   \ifcase \lastpenalty
        \then   \shortenproofright=0pt  
        \or     \unpenalty\hfil         
        \or     \unpenalty\unskip       
        \else   \shortenproofright=0pt  
        \fi
\fi
\global\dimen0=\shortenproofleft
\global\dimen1=\shortenproofright
\global\dimen2=\proofrulebreadth
\global\dimen3=\proofbelowshift
\global\dimen4=\proofdotseparation
\global\count255=\proofdotnumber
$\egroup  
\shortenproofleft=\dimen0
\shortenproofright=\dimen1
\proofrulebreadth=\dimen2
\proofbelowshift=\dimen3
\proofdotseparation=\dimen4
\proofdotnumber=\count255
}

\def\proofover{
\eproofbit 
\setbox\proofbelow=\hbox\bgroup 
\let\wereinproofbit\proofover
$\displaystyle
}%
\def\proofoverdbl{
\eproofbit 
\doubleprooftrue
\setbox\proofbelow=\hbox\bgroup 
\let\wereinproofbit\proofoverdbl
$\displaystyle
}%
\def\proofoverdots{
\eproofbit 
\proofdotstrue
\setbox\proofbelow=\hbox\bgroup 
\let\wereinproofbit\proofoverdots
$\displaystyle
}%
\def\proofusing{
\eproofbit 
\setbox\proofrulename=\hbox\bgroup 
\let\wereinproofbit\proofusing
\kern0.3em$
}

\def\endprooftree{
\eproofbit 
  \dimen5 =0pt
\dimen0=\wd\proofabove \advance\dimen0-\shortenproofleft
\advance\dimen0-\shortenproofright
\dimen1=.5\dimen0 \advance\dimen1-.5\wd\proofbelow
\dimen4=\dimen1
\advance\dimen1\proofbelowshift \advance\dimen4-\proofbelowshift
\ifdim  \dimen1<0pt
\then   \advance\shortenproofleft\dimen1
        \advance\dimen0-\dimen1
        \dimen1=0pt
        \ifdim  \shortenproofleft<0pt
        \then   \setbox\proofabove=\hbox{%
                        \kern-\shortenproofleft\unhbox\proofabove}%
                \shortenproofleft=0pt
        \fi
\fi
\ifdim  \dimen4<0pt
\then   \advance\shortenproofright\dimen4
        \advance\dimen0-\dimen4
        \dimen4=0pt
\fi
\ifdim  \shortenproofright<\wd\proofrulename
\then   \shortenproofright=\wd\proofrulename
\fi
\dimen2=\shortenproofleft \advance\dimen2 by\dimen1
\dimen3=\shortenproofright\advance\dimen3 by\dimen4
\ifproofdots
\then
        \dimen6=\shortenproofleft \advance\dimen6 .5\dimen0
        \setbox1=\vbox to\proofdotseparation{\vss\hbox{$\cdot$}\vss}%
        \setbox0=\hbox{%
                \advance\dimen6-.5\wd1
                \kern\dimen6
                $\vcenter to\proofdotnumber\proofdotseparation
                        {\leaders\box1\vfill}$%
                \unhbox\proofrulename}%
\else   \dimen6=\fontdimen22\the\textfont2 
        \dimen7=\dimen6
        \advance\dimen6by.5\proofrulebreadth
        \advance\dimen7by-.5\proofrulebreadth
        \setbox0=\hbox{%
                \kern\shortenproofleft
                \ifdoubleproof
                \then   \hbox to\dimen0{%
                        $\mathsurround0pt\mathord=\mkern-6mu%
                        \cleaders\hbox{$\mkern-2mu=\mkern-2mu$}\hfill
                        \mkern-6mu\mathord=$}%
                \else   \vrule height\dimen6 depth-\dimen7 width\dimen0
                \fi
                \unhbox\proofrulename}%
        \ht0=\dimen6 \dp0=-\dimen7
\fi
\let\doll\relax
\ifwasinsideprooftree
\then   \let\VBOX\vbox
\else   \ifmmode\else$\let\doll=$\fi
        \let\VBOX\vcenter
\fi
\VBOX   {\baselineskip\proofrulebaseline \lineskip.2ex
        \expandafter\lineskiplimit\ifproofdots0ex\else-0.6ex\fi
        \hbox   spread\dimen5   {\hfi\unhbox\proofabove\hfi}%
        \hbox{\box0}%
        \hbox   {\kern\dimen2 \box\proofbelow}}\doll%
\global\dimen2=\dimen2
\global\dimen3=\dimen3
\egroup 
\ifonleftofproofrule
\then   \shortenproofleft=\dimen2
\fi
\shortenproofright=\dimen3
\onleftofproofrulefalse
\ifinsideprooftree
\then   \hskip.5em plus 1fil \penalty2
\fi
}

\makeatletter
\def\section{\@startsection{section}{1}{0pt}{-3.25ex plus -1ex minus 
-.2ex}{1.5ex plus .2ex minus .3ex}{\normalfont\large\bf}}
\def\subsection{\@startsection {subsection}{2}{0pt}{-2ex plus -1ex minus 
   -.2ex}{1.5ex plus .2ex minus .3ex}{\@setfontsize\normalsize\@xpt{12}\bf}}
\makeatother

\newcommand{\defn}[1]{{\textit{\textbf{#1}}}}
\newtheorem{theorem}{Theorem}
\newcommand{\ie}{\emph{i.e.}}
\newcommand{\cf}{\emph{cf.}}

\title{%
\vspace{-7ex}
{\Large\bf Classical logic \;\;=\;\; Fibred MLL\vspace*{-.5ex}%
}}%
\author{\normalsize
   \sc Dominic Hughes \\ 
   \small Stanford University 
}
\date{\vspace{-1ex}\small March 25, 2005\thanks{Accepted for a short presentation at \emph{Logic in Computer Science '05}.}}
\begin{document}

\advance\textheight 1.15em
\maketitle
\thispagestyle{empty}
\vspace{-2ex}%
Syntactically, classical logic decomposes thus \cite{Gir87}:
\begin{center}
\sl Classical logic \;\;=\;\; MLL \;+\; Superposition
\end{center}
where MLL is Multiplicative Linear Logic, and superposition means
contraction (binary case) and weakening (nullary case).
Proof nets for classical logic have been proposed \cite{Gir91,Rob03},
containing explicit contraction and weakening nodes.  Just as the
logical sequent rules of conjuction (tensor) and disjunction (par)
are represented explicitly in the net, so are the structural rules,
contraction and weakening.
Unfortunately, this explicit representation of structural rules
retains some of the redundancy (or `syntactic bureaucracy') of sequent
calculus.  For example, in the net presentation of classical
categories in \cite{FP05} one has to quotient by equations on nets,
involving contraction and weakening.
Thus, proof nets for classical logic fail to achieve the same elegance
as the proof nets for MLL.

We present a representation of a proof which is more abstract than a
proof net, which we call a \emph{combinatorial proof}: superposition
is modelled \emph{mathematically}, as a lax form of fibration, rather
than syntactically (as in proof nets, which involve contraction and
weakening nodes).
This draws a nice boundary between logical rules and structural rules:
the former are modeled with explicit nodes, the latter as actual
superposition, in an abstract mathematical sense, without nodes.
We can summarise the situation thus:
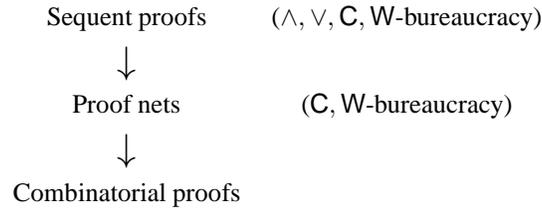
\vspace{-1ex}
\begin{center}
\begin{tabular}{cc}
Sequent proofs
& ($\wedge,\vee,\textsf{C},\textsf{W}$-bureaucracy)\\[1ex]
\Large$\downarrow$\\[.5ex]
Proof nets
& ($\textsf{C},\textsf{W}$-bureaucracy)\\[1ex]
\Large$\downarrow$\\[.5ex]
Combinatorial proofs
\end{tabular}
\end{center}
The arrows represent the abstraction from sequent calculus to proof
nets, and from proof nets to combinatorial proofs.  To the right, we
list the form of `bureaucracy' present at that layer, that is, the
connectives involved in representation redundancy, \ie, sensible rule
transpositions for identifying proofs.  ($\textsf{C}$ and $\textsf{W}$
stand for contraction and weakening.)

The definition of combinatorial proof is very simple, so we can give
it in full here.  Let $A$ be a formula of propositional classical
logic, generated by $\wedge$ and $\vee$ from literals (variables
$P,Q,\ldots$ and their negations $\overline{P},\overline{Q},\ldots$).
We identify a formula $A$ with its parse tree, a tree labelled with
literals at the leaves and $\wedge$ or $\vee$ at internal vertices. An
\defn{$\vee$-resolution} (or \defn{$\vee$-strategy}) 
of $A$ is any result of deleting one argument subtree from each $\vee$
(\cf\ \cite{HG03}).  A set $L$ of leaves of $A$ is a \defn{clique} of
$A$ if it is the set of leaves of a $\vee$-resolution of $A$.  A
\defn{combinatorial proof} of $A$ is an MLL proof net on a formula
$A'$ (viewed as an MLL formula, \ie, $\wedge$ = tensor, $\vee$ = par),
equipped with a label- and clique-preserving function from the leaves
of $A'$ to the leaves of $A$.  For example, here is a combinatorial proof of
Peirce's law, $((P\Rightarrow Q)\Rightarrow P)\Rightarrow
P\;=\;((\overline{P}\vee Q)\wedge
\overline{P})\vee P$\,:
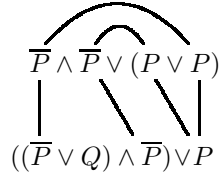
\begin{displaymath}\thicklines
\begin{array}{c}
\\[1ex]
\overline{P}%
        \begin{picture}(0,0)%
                \put(-2,11){\qbezier(0,0)(27,30)(54,0)}%
                \put(-4,-3){\line(0,-1){21}}
        \end{picture}%
\wedge \overline{P}%
        \begin{picture}(0,0)%
                \put(-2,11){\qbezier(0,0)(9,13)(18,0)}%
                \put(0,-3){\qbezier(0,0)(6,-10.25)(12,-20.5)}%
        \end{picture}%
\vee (P%
        \begin{picture}(0,0)%
                \put(-2,-3){\qbezier(0,0)(6,-10.25)(12,-20.5)}%
        \end{picture}%
\vee P%
        \begin{picture}(0,0)%
                \put(-4,-3){\line(0,-1){21}}
        \end{picture}%
) \\[5ex]
\!\!\!\!\!\!((\overline{P}\vee Q)\wedge\overline{P})\!\vee\! P
\end{array}
\end{displaymath}
The MLL proof net is drawn above, with two axiom links.  The label-
and clique-preserving function is shown by the
vertical lines.\vspace{-1ex}
\begin{theorem}[Soundness and Completeness]
A formula of classical propositional logic is true (valid)
\emph{iff} it has a combinatorial proof.
\end{theorem}
The label- and clique-preserving function is a lax form of fibration.
(We prove this via \cite{Hug04b}).  Thus a combinatorial proof is a
`fibred' MLL proof net, and we obtain the slogan in the title of the
paper:
\begin{center}
\sl Classical logic \;\;\;=\;\;\; Fibred MLL
\end{center}
Any sequent calculus proof translates to a combinatorial proof in a
relatively obvious way.  For a suitable variant of sequent calculus,
one obtains a surjection (\ie, sequentialisation theorem).  Cut
elimination for combinatorial proofs retains the richness of sequent
calculus: its non-determinism is not confluent.
In contrast, direct translation of a sequent proof to a linking (the
set of axiom links traced down from the axioms) leads to a certain
degeneracy: cut elimination becomes confluent \cite{LS05}.  Also in
contrast to \cite{LS05}, the MIX rule in this setting is optional.

Combinatorial proofs, in a somewhat different guise, were first
presented in \cite{Hug04b}, with MIX implicit.

\end{document}